\documentstyle[11pt]{article}
\begin{document}

\title{\textbf{Rank-one Solutions for Homogeneous Linear Matrix Equations over the
Positive Semidefinite Cone  } }

\author{Yun-Bin Zhao \thanks{School of
Mathematics, University of Birmingham, Edgbaston B15 2TT,
Birmingham,  United Kingdom ({\tt y.zhao.2@bham.ac.uk}). } ~ and Masao
Fukushima \thanks{Department of Applied Mathematics and Physics,
Graduate School of Informatics, Kyoto University, Kyoto 606-8501,
Japan ({\tt fuku@i.kyoto-u.ac.jp}).}}

\date{ (Revised, September 2012) }

\maketitle

{\bf Abstract. }   The problem of finding a rank-one
 solution to a system of linear matrix
equations arises from many practical applications.  Given a system
of linear matrix equations, however, such a low-rank solution does
not always exist. In this paper, we aim at developing some
sufficient conditions for the existence of a rank-one solution to
the  system of homogeneous linear matrix equations  (HLME) over the
positive semidefinite cone. First,  we prove that an existence
condition of a rank-one solution can be established by a homotopy
invariance theorem. The derived condition is closely related to the
so-called $P_\emptyset$ property of the function defined by
quadratic transformations. Second, we prove that the existence
condition for a rank-one solution  can be also established through
the maximum rank of the (positive semidefinite) linear combination
of given matrices. It is shown that an upper bound for the rank of
the solution to a system of HLME  over the positive semidefinite
cone can  be obtained efficiently by solving a semidefinite
programming (SDP) problem. Moreover, a sufficient condition for the
nonexistence of  a rank-one solution to the system of HLME is also
established in this paper. \\


 {\bf AMS. } 90C22, 90C25, 65K05, 15A60.

\newpage

\section{Introduction} Let $R^n$ denote the
$n$-dimensional Euclidean space, with the standard inner product,
and let $S^n$ denote the set of real symmetric matrices. For a given
$A\in S^n$, $A\succeq 0  ~(\succ 0)$ means that $A$ is positive
semidefinite (positive definite). For two $n\times n$ matrices $X$ and $Y$,
$\langle X, Y \rangle =\textrm{tr}(X^T Y) $ denotes the
inner product of $X$ and $Y$, where $\textrm{tr}(\cdot)$ stands for the
trace of a square matrix. We use $\|X\|$ and $\|X\|_*$ to denote the spectral norm
and the nuclear norm (i.e., the sum of singular values), respectively, of matrix $X$.

Stimulated by the recent work on compressed/compressive sensing
(e.g. \cite{DH01, D06, CRT06, C06}), the study of finding a low-rank
solution to optimization problems with linear matrix (in)equality
constraints has recently become intensive  \cite{RFP07, Z12LAA}.
Many practical problems across disciplines (such as sparse signal
recovery \cite{D06,C06}, system control \cite{EGG93, EG94, HTS99,
F02, FHB04}, matrix completion \cite{CR08, CT09}, machine learning
\cite{ABEV06, MJCD08}, quadratic equation \cite{Z12LAA}, Euclidean distance geometry \cite{T00, SY07,
D09} and combinatorial optimization \cite{AV09}) can be formulated
as the following problem:
$$\min  \left\{ \textrm{rank}(X):  ~
{\mathcal A}  X =b,  ~  X \in  C  \right\},$$ where $X$ is an $n
\times p$ matrix, $ b$ is a vector in $R^m,$ $ {\mathcal A}:  R^{n
\times p} \to R^m$ is a linear operator, and $C \subseteq R^{n
\times p}$ is a convex set. Various heuristic methods for such
problems have been proposed and investigated (e.g. \cite{F02, RFP07,
MGC09, TY09, AI10}). Among all low-rank solutions, a rank-one
solution is particularly useful in many situations, especially in
system control and quadratic optimization \cite{EGG93, EG94, HTS99,
AV09, H01}. Locating a low-rank solution (especially a rank-one
solution) is not only motivated by these practical applications, but
also motivated naturally by the structure of the solution set of
semidefinite programming (SDP) problems. In fact, an optimal
solution of an SDP problem, if attained, usually lies on the
boundary of its feasible set. Note that the nontrivial extreme rays
of the positive semidefinite cone are generated by rank-one
matrices.

However, for a given system of linear matrix equalities, a rank-one
solution with a desired structure does not always exist. In this
paper, we focus on the existence issue of a rank-one solution to the
homogeneous system:
\begin{equation}\label{LME-PSDC} \langle A_i, ~X\rangle =0, ~
i=1,\dots,m , ~ X\succeq 0, \end{equation}  where $A_i \in S^n, i=1,
\dots,m, $ are given $n\times n$ matrices. This system is referred
to as the \emph{homogeneous linear matrix equations (HLME) over the
positive semidefinite cone}.  We aim at addressing the following
fundamental question:   When does system  (\ref{LME-PSDC}) possess a
rank-one solution?  The main motivation to the study of the system
(\ref{LME-PSDC}) is that it is closely related to the system of
quadratic equations (see e.g., \cite{Z12LAA}):
\begin{equation}
\label{quadratic}  x^T A_i x=0, ~i=1, \dots,m, ~x\in R^n.
\end{equation}
Hence it is closely linked to various linear algebra and
optimization topics, such as the simultaneous diagonalizability of a
given set of matrices $(A_1,\dots, A_m)$,  the convexity property of
the field of values associated with a finite number of matrices
\cite{B61, HJ91, HUT02, HU07}, and the widely used `S-Lemma' or
`S-Procedure' in system control and optimization
\cite{BV04,EG94,BN01,PT07, P98}.

  In this paper, we investigate the
existence of a rank-one solution to the system (\ref{LME-PSDC}) from
the viewpoint of nonlinear analysis and rank optimization. Clearly,
the system (\ref{LME-PSDC}) has a rank-one solution if and only if
(\ref{quadratic}) has a nonzero solution.
 Thus, the existence of a rank-one solution of the system (\ref{LME-PSDC}) is
equivalent to that of a nonzero solution of the system
(\ref{quadratic}). Except for some special cases, however, a general
and complete characterization of the existence condition
 for a rank-one solution to the system (\ref{LME-PSDC})
remains open (see, for instance, the open ``Problem 12" and
``Problem 13" in \cite{HU07}). The study of (\ref{quadratic}) can
date back to 1930s. Thanks to the early work of Dines \cite{D41,
D42, D43}, Brickman \cite{B61}, Calabi \cite{C64}, and the work of
Finsler \cite{F37},   a complete characterization of the system
(\ref{quadratic}) with $m=2$ and $n\geq 3$ is clear: \emph{$x=0$ is
the only solution to the system $x^T A_1 x=0$ together with $ x^T
A_2 x =0 $  if and only if $t_ 1 A_1+t_2 A_2 \succ 0 $ for some
$t_1, t_2\in R.$}   A good survey on the historical development of
this result can be found in \cite{U79, HUT02, P98}. This result is
closely related to S-Lemma/S-procedure \cite{PT07}, and related to
the approximate S-lemma and low rank issues discussed in
\cite{Bar02}. It can
  be related to the trust region subproblem in nonlinear
optimization as well (see e.g. \cite{SW95, NW99}).

Let us restate the above classical result as follows, in terms of
rank-one solutions to the system of HLME over the positive
semidefinite cone. We  call it ``Dines-Brickman's Theorem".\\

\textbf{Theorem 1.1 } (Dines-Brickman)  \emph{When $m=2$ and $n\geq 3,$
the system
\begin{equation}\label{AABB} \langle A_1, X\rangle =0, ~\langle A_2,
X\rangle =0, ~ X\succeq 0 \end{equation}  has a rank-one solution if
and only if $t_1 A_1+t_2 A_2 \not \succ 0 $ for any $t_1, t_2
\in R,$  in other words,
$$\max_{t_1, t_2}~\{ \, {\rm rank} (t_1A_1+t_2 A_2):
~t_1A_1+t_2 A_2 \succeq 0\} \, \leq \, n-1.$$}

 Unfortunately, such a
complete characterization does not hold in general for $m\geq 3.$
The first purpose of this paper is to establish a general sufficient
condition for the existence of a rank-one solution to the system
(\ref{LME-PSDC}) by using a homotopy invariance theorem. To this
end, we introduce a class of functions called $P_\emptyset,$ and
show that the system (\ref{LME-PSDC}) has a rank-one solution if the
quadratic image function has a $P_\emptyset$ property.  This
analysis is the first to link $P_\emptyset$-functions and the
existence of a rank-one solution   to the system of HLME over the
positive semidefinite cone. The second purpose of this paper is to
develop an existence condition  for the rank-one solution of
(\ref{LME-PSDC}) via a rank optimization approach. To this end, we
introduce  the following rank maximization problem
$$
r^*= \max_{t_1,\dots,t_m\in R} \left\{ \textrm{ rank }
\left(\sum_{i=1}^m t_i A_i \right) \, : ~~ \sum_{i=1}^m t_i
A_i\succeq 0 \right\}, $$ which turns out to be an important factor
for the existence of a rank-one solution of the system
(\ref{LME-PSDC}) (this has already been observed in the
aforementioned Theorem 1.1).
  In section 3, we also point out that the value of $r^*$, combined with
    Barvinok-Pataki's bound \cite{B95, PA98}, can be used to determine
an upper bound for the rank of the solution of (\ref{LME-PSDC}).

It is also important to understand when the system (\ref{LME-PSDC})
does not have a rank-one solution. In other words, we study the
question: When is $x=0$ the only solution to the quadratic system
(\ref{quadratic})? This was posted as the open ``Problem 13" in
\cite{HU07}. For $m\geq 3$, the well-known condition $\sum_{i=1}^m
t_i A_i \succ 0 $ does imply that $x=0$ is the only solution of the
quadratic system (\ref{quadratic}) (and hence the system
(\ref{LME-PSDC}) has no rank-one solution). However, this condition
is too strong (see Lemma 3.2 for details). This motivates us to
investigate  the above-mentioned open question, and to develop
another sufficient condition for the nonexistence of a rank-one
solution to the system (\ref{LME-PSDC}) from a nonlinear analysis
point of view.

 This paper is organized as follows. In section~2, we develop
a sufficient condition for the existence of a rank-one solution to
(\ref{LME-PSDC}) by using a homotopy invariance theorem. In
section~3, we establish other sufficient conditions by means of rank
optimization. A nonexistence condition for the rank-one solution of
(\ref{LME-PSDC}) is provided in section~4, and conclusions are given
in the last section.

 \section{Existence of a rank-one solution: homotopy invariance}

 Note that any rank-one matrix in a positive semidefinite cone must be
of the form  $X=x x^T$ with $x\not=0.$ Thus the system
(\ref{LME-PSDC}) has a rank-one solution if and only if there is an
$x\not=0$ such that $\langle A_i, xx^T\rangle=0, \, i=1,\dots, m,$
which is nothing but the system (\ref{quadratic}). By homogeneity,
this is equivalent to saying that the system (\ref{LME-PSDC}) has a
rank-one solution if and only if there is a solution to the system
\begin{equation} \label{1111}  x^TA_i x
=0, ~i=1,\dots, m, ~ x^Tx=1, \end{equation}
which has $m+1$ equations and $n$
variables. It is worth mentioning that the system (\ref{1111}) is
 related to the structure of the set $\{x: ~x^T A_ix \leq 1, ~i=1,\dots, m\}$
 which has been considered by Nemirovski, Roos and Terlaky \cite{NRT99}.
 For instance, if (\ref{1111}) has no solution, then
the set $\{x: ~x^T A_ix \leq 1, ~i=1,\dots, m\}$ contains only
a single point $x=0.$

 In this section, we assume that $m+1\leq n.$ Let $\Phi:
R^n \to R^{m+1} $ be the mapping defined
 by
 \begin{equation} \label{Phi} \Phi(x) = (x^T A_1 x, ~
     \cdots, ~
 x^T A_m x, ~x^T x-1 )^T. \end{equation}
 When $m+1<n $, we may introduce the extra matrices $
A_{m+1}= \cdots = A_{n-1}= 0$ into (\ref{LME-PSDC}) without any
change of the system. Thus, without loss of generality, we assume
that $ m+1 =n $ in the remainder of this section
so that $\Phi(x)=0$ consists of $n$ equations with $n$ variables.
  An immediate observation is given as follows.

  \vskip 0.08in

\textbf{Lemma 2.1.}  \emph{The system (\ref{LME-PSDC}) has a
rank-one solution if and only if $\Phi(x)=0$ has a solution.}

  \vskip 0.08in

This observation combined with the next result (homotopy invariance
theorem of topological degree) is a key to developing our first
existence result. For a bounded subset $ D $ of $ R^n, $
  $\overline{D}$ and $\partial D $ denote
 the closure  and the boundary of $D,$ respectively. Let $f$ be a continuous  function from
 $\overline{D} $ into $ R^n$. For $ y\in R^n$  such that $ y\not\in f(\partial D)=\{f(x): x\in \partial D \},$
 deg$(f,D, y) $ denotes the topological degree associated with $ f,D$ and $ y,
 $ which has been widely used in the existence analysis  of
 nonlinear equations
(see, \cite {OR70, L78, OCC06}).  The following result is called the
 homotopy invariance theorem. Part (i) is due to Poincar\'e and Bohl,
while part (ii) is attributed to Kronecker \cite[Theorems 6.2.4 and
6.3.1]{OR70}.

  \vskip 0.08in

\textbf{ Lemma 2.2.} (Poincar\'e-Bohl-Kronecker) \emph{ Let $ D\subset R^n $
be a nonempty open bounded set and $ F, G $ be two continuous functions from
$\overline {D}$ into $R^n$. Let  $H(x,t)=tG(x)+(1-t)F(x),$ $ 0\leq
t\leq 1,$ and let $ y $ be an arbitrary point in $ R^n $. }

(i) \emph{If $ y\notin\{H(x,t): x\in \partial D \mbox{ and } t\in [0,1]\},$
then $ \mbox{\rm deg}(G, D,y)=\mbox{\rm deg}(F,D,y).$}

(ii)  \emph{If $ y\notin F(\partial D) $ and  \mbox{\rm deg}$(F, D,
y)\not= 0$, then the equation  $ F(x)=y $ has a solution in $ D$.}

  \vskip 0.08in

First, we prove the following technical result.

  \vskip 0.08in

\textbf{Lemma 2.3.} \emph{If $\Phi(x)=0$ has no solution, where $\Phi$ is
defined by (\ref{Phi}), then there exists a sequence $\{x^k\} \subset R^n$
such that $\|x^k\|\to \infty $ as $k\to \infty ,$ and
for each $k$ there exists $\gamma^k\in (0,\infty)$ such that $ \Phi
(x^k) =-\gamma^k x^k, $ i.e., $x^k $ satisfies the following
relations:
\begin{eqnarray}
  (x^k)^T A_i x^k  & = &  - \gamma^k x^k_i,   ~~i=1,\dots, m, \label{TH111} \\
  (x^k)^T x^k & = & 1-\gamma^k x^k_{n}, \label{TH222}
  \end{eqnarray}
where $x^k_j$ denotes the $j$th component of $x^k.$}

  \vskip 0.08in

 \emph{Proof.} Suppose that $\Phi(x)=0$ has no solution.
Consider the homotopy between the identity
mapping and $\Phi (x) $, i.e.,
$$ {\cal H}  (x,t)= t x +(1-t)\Phi (x) = \left(\begin{array}{c}
t x_1 +(1-t)  x^T A_1 x    \\
     \vdots  \\
tx_m + (1-t) x^T A_m x \\
t x_{n}+ (1-t) (x^Tx -1)   \\
\end{array}
\right)\in R^{n}.
$$
Let
$$ {\mathcal Q}  = \{x \in R^{n}: ~ {\cal H}(x,t)=0\textrm{  for some  } t \in [0,1] \} .  $$
First, we prove that  $\mathcal{Q} $ is unbounded. In fact, if
${\mathcal Q}$ is bounded, then there exists an open bounded ball $D \subset R^n$
that contains the set $\mathcal{Q} , $ and
$D$ can be chosen large enough so that the boundary  $\partial D$
does not touch the set ${\mathcal Q}, $ i.e., $\partial D\bigcap
\mathcal{Q} =\emptyset.$ By the definition of $\mathcal{Q} $, we
deduce that
$$0 \notin\{{\cal H}(x,t): ~~x\in \partial D, ~~0\leq t\leq 1\}. $$
(Indeed, if $0$ is an element of the above set, then ${\cal
H}(x,t)=0$ for some $x\in\partial D $ and $t\in [0,1],$ which
implies $\partial D\bigcap \mathcal{Q}  \not= \emptyset,$ a
contradiction.) Thus, $\deg (I, D,0) $ and $\deg(\Phi, D,
0)$ are well defined, and by Lemma 2.2 (i), we have
$$
\deg(\Phi, D, 0)=  \deg (I, D,0). $$ For the identity mapping, we
have $|\deg (I, D,0)|=1.$ Thus, we have  $ \deg(\Phi, D, 0) \not =
0$, which along with Lemma 2.2 implies that there exists a solution to $
\Phi (x)=0,  $ contradicting the assumption. So $\mathcal{Q} $
must be unbounded, and  there exists an unbounded sequence $\{x^k\}$
in ${\mathcal Q}.$   Without loss of generality, let $x^k \ne 0$ for
all $k,$ and $ \|x^k\| \to \infty $ as $k\to \infty. $ Since $
\{x^k\} \subseteq {\mathcal Q}$,   there is a sequence
$\{t^k\}\subseteq [0,1] $ such that
\begin{equation}\label{HHH} {\cal H} (x^k, t^k) = \left(\begin{array}{c}
t^k x^k_1 +(1-t^k)  (x^k) ^T A_1 x^k   \\
     \vdots  \\
t^kx^k_{m} + (1-t^k)  (x^k)^T A_m x^k\\
t^k
x^k_{n}+ (1-t^k)  ((x^k)^Tx^k-1)  \\
\end{array}
\right) =    0 .  \end{equation}
  Since $x^k \ne 0,$   it follows from (\ref{HHH})  that $ t^k\not=1.$  By assumption,
   there is no solution to the equation
  $\Phi(x) =0,$ so it follows
from (\ref{HHH}) that
  $t^k\not=0.  $ As a result,  (\ref{HHH}) can be written as
$$ \Phi  (x^k) = -\left(\frac{t^k}{1-t^k}\right) x^k, ~~ t^k \in (0,1), ~~ k \geq 1.  $$
By the definition of $\Phi,$ we have
\begin{eqnarray*}
(x^k)^T A_i x^k &  = & -\frac{t^k} { 1-t^k} x^k_{i} , ~~ i=1, \dots, m, \\
(x^k)^T x^k  & = & 1 - \left(\frac{t^k }{1-t^k} \right) x^k_{n},
\end{eqnarray*}
where $t_k\in (0,1) $ for all  $ k\geq 1. $  The desired result then
follows by setting $\gamma^k= \frac{t^k }{1-t^k}. $ ~~ $  \Box$ \\

A similar analysis to the above (by degree theory) has been used in
the existence analysis for the solution of finite-dimensional
variational inequalities and complementarity problems (see e.g.
\cite{IZ00a, ZH99, ZI00,ZL01}). It is also worth mentioning that the
identity mapping in ${\cal H}(x,t)$ can be replaced by a general
invertible mapping $ \varphi(x) $ with $|\textrm{deg} (\varphi, D,
0)|\not=0, $  and thus the result developed in this paper can be
easily adapted to this case. However, we choose to use the identity
mapping throughout this paper  in order to keep  results as simple
as possible.

  \vskip 0.08in

From Lemma 2.3, we can prove the next result.

  \vskip 0.08in

\textbf{Lemma 2.4.} \emph{If $\Phi(x)=0$ has no solution, where $\Phi$ is
defined by (\ref{Phi}), then there exists a  vector $\widehat{x}$
with $\|\widehat{x}\|=1$ such that
\begin{equation}\label{accumulation-01} \widehat{x}^T A_i \widehat{x}=-\frac{1}{\delta} \widehat{x}_i, ~~ i=1,\dots,
m, ~  -\widehat{x}_n= \delta \in (0,1),\end{equation}
 and hence
\begin{equation}\label{accumulation-02} \widehat{x}^T \left( \sum_{i=1}^m  -\widehat{x}_i A_i\right)  \widehat{x}
=\frac{1}{\delta} (1-\delta^2) >0.
\end{equation}}

\emph{Proof.} Suppose that $\Phi(x)=0$ has no solution. Then, by
Lemma 2.3, there
 exists a sequence $\{x^k\}$ satisfying (\ref{TH111}) and
(\ref{TH222}),  and $\|x^k\|\to \infty $ as $k\to \infty.$  We see
from (\ref{TH222}) that $ x^k_n <0  $ for all sufficiently large
$k.$  It is not difficult to show that
 any accumulation point of the sequence $ \left\{  x^k_n/\|x^k\| \right\} $
is in $(-1,0). $ In fact, since  $x^k_n$ is negative for all
sufficiently large $k,$ any accumulation point of the sequence $
\left\{ x^k_n/ \|x^k\| \right\}$ must be in $[-1, 0].$  So it is
sufficient to prove that it is not equal to $ -1$ or $0.$  Let
$\widehat{x}$ be an arbitrary accumulation point of the sequence $
\left\{ x^k/\|x^k\| \right\}$. By passing to a subsequence if
necessary, we may assume that $\lim_{k\to \infty} x^k /\|x^k\|
=\widehat{x}.$ By (\ref{TH222}), we have
\begin{equation}\label{111}
 1= \frac{1}{\|x^k\|^2}-\frac{\gamma^k}{\|x^k\|} \left(\frac{x^k_n
}{\|x^k\|}\right)
\end{equation}
and, by (\ref{TH111}), for every $i=1, \dots,m $, we have
\begin{equation}\label{rrr}
\left(-\frac{\gamma^k}{\|x^k\|}
\right)\left(\frac{x^k_i}{\|x^k\|}\right)=\frac{(x^k)^T A_i
x^k}{\|x^k\|^2}\to \widehat{x}^T A_i \widehat{x}.
\end{equation}

\emph{Case 1:} Assume that $\widehat{x}_n = -1. $  Then,  from (\ref{111}), we see that $
\gamma^k/\|x^k\| \to 1.$ Since $\widehat{x}_n = -1$ and
$\|\widehat{x}\| =1,$ we must have $ x^k_i/\|x^k\|\to \widehat{x}_i
=0 $ for all $ i=1, \dots, m.$ Thus, it follows from (\ref{TH111})
that
 $$ \widehat{x}^TA_i\widehat{x} = -\widehat{x}_i =0,  ~~ i=1, \dots, m. $$
This implies that the unit vector $\widehat{x}$ is a  solution of
the system
  (\ref{quadratic}), i.e., $\Phi(\widehat{x})=0, $  contradicting the
assumption of this lemma.

\emph{Case 2:} Assume that $\widehat{x}_n = 0. $
Then it follows  from (\ref{111}) that $
 \gamma^k/\|x^k\| \to \infty $ as $ k\to\infty $.
Thus, from (\ref{rrr}) we  have  $ x^k_i/\|x^k\| \to \widehat{x}_i
=0,$  $  i=1, \dots, m. $ This contradicts the fact that $\|
\widehat{x} \|=1. $

Both    cases above yield a contradiction. Thus we conclude that for
any accumulation point $\widehat{x}$ of $\{x^k /\|x^k\|\}$, its last
component $\widehat{x}_n $ satisfies $\widehat{x}_n \in (-1,0).$ We
now show that such an accumulation point  satisfies
(\ref{accumulation-01}) and (\ref{accumulation-02}). Indeed, without
loss of generality,
  we assume that
 $ \lim_{k\to \infty} x^k/\|x^k\| =  \widehat{x}$ and
$- \widehat{x}_n = \delta \in (0,1). $ From (\ref{111}), we see that
$ \gamma^k/\|x^k\| \to 1/\delta,$   and thus (\ref{accumulation-01})
follows directly from (\ref{rrr}). Multiplying  (\ref{TH111}) by
$x^k_i$ for each $i$ and adding them up yield
$$ (x^k)^T \left(\sum_{i=1}^m x^k_i A_i \right) x^k =-\gamma^k
\sum_{i=1}^m (x^k_i)^2 =-\gamma^k \left(\|x^k\|^2-(x^k_n)^2\right).
$$ Dividing this equality by $\|x^k\|^3$ and taking the limit yield
(\ref{accumulation-02}). ~~  $ \Box $

\vskip 0.08in

An immediate consequence of Lemmas 2.1 and 2.4 is the following
sufficient condition.

  \vskip 0.08in

\textbf{Corollary 2.5.} \emph{If the system
$$ x^T A_i x= \frac{  x_i}{x_n}, ~~ i=1,\dots, m,
~ x_n  \in (-1,0), ~\|x \| =1 $$ is inconsistent, then the system
(\ref{LME-PSDC}) has a rank-one solution.}

\vskip 0.08in

Based on this fact, we have the following result.

\vskip 0.08in

 \textbf{Theorem 2.6.} \emph{If $A_i, i=1, \dots,m, $ satisfy the condition
\begin{equation}\label{P0-condition}
  \max_{1\leq i\leq m, ~x^T A_i x \not =0}  x_i  (x^T A_i x) \geq
 0 ~\textrm{ for all  } x\textrm{ with }\|x\| =1\textrm{ and } x_n<0,
 \end{equation}
then the system (\ref{LME-PSDC}) has a rank-one solution.}

  \vskip 0.08in

\emph{ Proof.} We prove this result by contradiction. Assume that
the system (\ref{LME-PSDC}) has no rank-one solution. Then, by Lemma
2.1, the equation $\Phi(x)=0 $  has no solution.  It follows from
Lemma 2.4 that there exists a unit vector   $\widehat{x}$ satisfying
(\ref{accumulation-01}).
 Multiplying both sides of
(\ref{accumulation-01}) by $\widehat{x}_i$ yields
$$ \widehat{x}_i (\widehat{x}^T A_i \widehat{x}) = - \frac{1}{\delta}
\widehat{x}_i^2, ~~ i=1, \dots, m , ~ -\widehat{x}_n= \delta\in (0,1).
$$ Note that $\|\widehat{x}\|=1$ and  $|\widehat{x}_n|=\delta <1. $
Thus $(\widehat{x}_1, \dots, \widehat{x}_m)\not =0, $ which implies that
$$ \max_{1\leq i\leq m,  ~\widehat{x}^T  A_i  \widehat{x} \not=0}
\widehat{x}_i (\widehat{x}^T A_i \widehat{x})  = \max_{1\leq i\leq
m,  ~ \widehat{x}_i \not=0} (- \widehat{x}_i^2/\delta ) <0. $$ This
contradicts (\ref{P0-condition}).  ~~   $  \Box $

\vskip 0.08in

\textbf{Corollary 2.7.} \emph{If $A_i, i=1, \dots,m, $ satisfy the condition
\begin{equation}\label{geometric}
\max_{1\leq i\leq m, ~x^T A_i x \not =0}  x_i  (x^T A_i x) \geq
 0 ~ \textrm{ for all }x\not=  0,
 \end{equation}
then the system (\ref{LME-PSDC}) has a rank-one solution.}

\vskip 0.08in

 Clearly, the condition (\ref{geometric}) is stronger than  (\ref{P0-condition}).
So it implies the existence of a rank-one solution to the system
(\ref{LME-PSDC}). Motivated by conditions (\ref{P0-condition}) and
(\ref{geometric}), we introduce the class of $P_\emptyset$ functions
defined as follows.

\vskip 0.08in

\textbf{Definition 2.8. } \emph{Let $ D \subseteq R^n $ and
$\widehat{x}\in D.$ A mapping $F: R^n \to R^n $ is said to be a
$P_\emptyset$-function at $\widehat{x}$ over  $D$ if $$\max_{1\leq
i\leq n, ~ F_i(x) \ne F_i(\hat{x})  } (x_i-\widehat{x}_i)
(F_i(x)-F_i(\widehat{x})) \geq 0 ~ \textrm{ for all } \, x\in D, \, x\not=
\widehat{x}.
$$}

\vskip 0.08in

Recall that for a given set $ D \subseteq R^n $ and $\widehat{x}\in
D, $ the mapping $G: R^n \to R^n $ is said to be a $P_0$-function at
$\widehat{x}$ over  $D$ if $\max_{1\leq i\leq n, x_i\not=
\widehat{x}_i} (x_i-\widehat{x}_i) (G_i(x)-G_i(\widehat{x})) \geq 0
~ \textrm{ for all } x\not= \widehat{x}$ and $ x\in D. $ The class
of $P_0$-functions has been widely used in  nonlinear analysis and
optimization (see e.g. \cite{PS92, FP03, ZI00, ZL01}). There are
some relationships between $P_\emptyset$-functions and
$P_0$-functions. In fact, when the inverse of a
$P_\emptyset$-function exists, it is easy to see that the inverse is
a $P_0$-function. By Definition 2.8, we can state the following
result.

\vskip 0.08in

\textbf{Theorem 2.9.}  \emph{Let  ${\mathcal F}: R^n \to R^n$ be
defined by $$ {\mathcal F} (x)  = \left( x^T A_1 x, ~\dots, ~x^T A_m
x, ~x^T x\right)^T $$ and let $U=\{x: \|x\| \leq 1 \}.$  If ${\cal
F}$ is a $P_\emptyset$-function at $x=0$ over $U$, then the system
(\ref{LME-PSDC}) has a rank-one solution.}

\vskip 0.08in

\emph{Proof.}   Since ${\cal F} $ is a $P_\emptyset $-function at
$x=0$ over $U,$  we have for any $x \in U \setminus \{0\}$ that
\begin{equation}\label{FFFF} \max_{1\leq i\leq n, ~ {\cal
F}_i(x)\ne{\cal F}_i(0)} x_i ({\cal F}_i(x)-{\cal F}_i(0)) \geq 0.
\end{equation}  In particular, for any $x$ such that $ \|x\|=1 $
 and $x_n<0$, we have $x_n ({\cal F}_n(x)-{\cal F}_n(0)) =x_n x^T x =x_n<0.$
So  (\ref{FFFF}) is reduced to
$$ \max_{1\leq i\leq m, ~x^T A_ix  \not= 0 }   x_i (x^T A_i x) \geq 0. $$
 By Theorem 2.6 or its corollary,
the system (\ref{LME-PSDC}) must have a rank-one solution. ~ $ \Box
$

\vskip 0.08in

 \textbf{Remark 2.10.} In this section, we have seen
that the degree-based analysis  can provide
 a sufficient condition for the system (\ref{LME-PSDC}) to have a rank-one
 solution. In section 5, we will further show  that such a sufficient condition
 is
 almost necessary for the system (\ref{LME-PSDC}) to have a
 rank-one solution. In section 3, we show that checking whether or not the
 system (\ref{LME-PSDC}) has a rank-one solution is equivalent to
 solving an SDP problem with a rank constraint, which is  clearly not an
 easy problem due to the rank constraint. It is well known that a
 general rank minimization problem are NP-hard \cite{RFP07, Z12LAA} since it includes
 the so-called cardinality minimization problem as a special case.
So roughly speaking, the level of difficulty for checking the
conditions developed in this paper, such as the ones in Corollary
2.5, Theorem 2.6 and Theorem 4.5,  are almost
 equivalent to that of the original system  (\ref{LME-PSDC}) which,
  except for the case of $m=2$ and $n\geq 3,$   is difficult in general.
   However, these conditions provide a new angle (from degree theory) to understand the
 system (\ref{LME-PSDC}). More
 interestingly, some verifiable sufficient conditions (from a rank optimization point of view)
  can be also developed
  for the system (\ref{LME-PSDC}) to have a rank-one solution, as
  shown in the next section.

\section{Existence of a rank-one solution: rank optimization}

We now study the existence of a rank-one solution to the system
(\ref{LME-PSDC}) from a rank optimization point of view. When $m=2$
and $n\geq 3$, Theorem 1.1 claims that the condition $t_1 A_1+t_2
A_2 \not\succ 0 $ is a complete characterization of the existence of
a rank-one solution to the system (\ref{LME-PSDC}). However, this
result does not hold when $m\geq 3.$  It is interesting to note that
a complex counterpart of such a result was given in \cite{PT07}. In
this section, we show that some conditions stronger than
$\sum_{i=1}^m t_i A_i \not\succ 0 $ are needed in order to ensure
the existence of rank-one solutions. Before we proceed, let us first
reformulate the problem as a rank constrained optimization problem.

 Consider the  following rank-constrained optimization problem:
 \begin{equation}\label{rank-con}
\max \{\langle I, X \rangle: ~~ \langle A_i, X \rangle = 0, ~i=1,
\dots, m, ~\textrm{rank}(X) \leq 1,  ~ X\succeq 0 \}. \end{equation}
If $X=0$ is the only feasible point to the problem, then the optimal
value of the  problem is $0.$ Otherwise, it is $\infty.$ Since the
system (\ref{LME-PSDC}) is homogeneous, normalizing the system does
not change its solvability  and the rank of its solutions. So adding
the constraint ${\rm tr}(X)\leq 1 $ (i.e., $\langle I, X \rangle
\leq 1$) to  (\ref{rank-con}) yields the following problem:
\begin{eqnarray}  \label {rank-prob}  z^*= \, & \max & \langle I, X \rangle \nonumber \\
& \textrm{s.t.} & \langle A_i, X \rangle = 0, ~i=1, \dots, m,  \nonumber\\
 & & \textrm{rank} (X) \leq 1,   \\
 && \textrm{tr}(X) \leq 1,  \nonumber\\
 & & X\succeq 0. \nonumber
 \end{eqnarray}
Since the optimal value $z^*$ is either 0 or 1, we have the
following observation. \\

\textbf{Lemma 3.1.} \emph{The system (\ref{LME-PSDC}) has a rank-one
solution if and only if $z^*=1 $ is the optimal value of
(\ref{rank-prob}).} \\

In other words, $z^*=0$ is the optimal value of (\ref{rank-prob}) if
and only if (\ref{LME-PSDC}) has no rank-one solution. Lemma 3.1
indicates that checking the existence of a rank-one solution to the
system (\ref{LME-PSDC}) is equivalent to solving an SDP problem with
the rank constraint $\textrm{rank}(X)\leq 1$, which is hard to solve
in general because of the discontinuity and  nonconvexity of
$\textrm{rank}(X). $  Based on the problem obtained
 by dropping the constraint
``$\textrm{rank}(X) \leq 1" $ from (\ref{rank-prob}), we have the
following result.  \\

\textbf{Lemma 3.2. } \emph{There exist  $t_i \in R, i=1,\dots,m,  $
such that $\sum_{i=1}^m t_i A_i \succ 0 $ if and only if $X=0$ is
the only solution  to the system (\ref{LME-PSDC}).  In other words,
$\sum_{i=1}^m t_i A_i \not\succ 0 $ for all $t_i\in R, i=1,\dots,m$,
if and only if the system (\ref{LME-PSDC}) has a nontrivial
solution,
i.e., a solution $X$ with ${\rm rank}(X)\geq 1. $} \\

\emph{Proof.}  The standard SDP duality theory \cite{H01} (or the
result in \cite{BW81})  can yield the result of this lemma. In fact,
let us consider the SDP problem
\begin{equation}  \label {SDP-00}   \max   \{\langle I, X \rangle:
~\langle A_i, X \rangle = 0, ~i=1, \dots, m ,~
  \textrm{tr}(X) \leq 1,~
  X\succeq 0\},
 \end{equation}
and its  dual problem
\begin{equation}  \label {DSDP-00}  \min  \left\{\alpha : ~
 \sum_{i=1}^m t_i A_i  +\alpha I  \succeq I,
 ~ \alpha \geq 0, \ t_i \in R, i=1,\dots,m \right\}.
 \end{equation}
Clearly, (\ref{DSDP-00}) satisfies the
  Slater's condition (for instance, $(\alpha,t_1,\dots,t_m)=
(2,0,\dots,0) $ is a strictly feasible point). The optimal value of
(\ref{DSDP-00}) is obviously finite. By the duality theory of
 semidefinite programming, both problems (\ref{SDP-00})  and
(\ref{DSDP-00}) have finite optimal values and  there is no duality
gap between them (i.e., their optimal values are equal). If there
exist $t_i, i=1,\dots,m $ such that
 $ \sum_{i=1}^m t_i A_i \succ 0 ,$ then  $\sum_{i=1}^m (\beta t_i) A_i \succeq I $ for some $\beta>0,$
 which means that
 the optimal value of the
dual problem (\ref{DSDP-00}) is $0. $  Thus the optimal value of
 (\ref{SDP-00}) is also $0, $ implying that $X=0$ is the only point satisfying the system
(\ref{LME-PSDC}). Conversely, if  $X=0$ is the only solution of
(\ref{LME-PSDC}), then the optimal value of (\ref{SDP-00}) is $0$,
and thus  the dual optimal value is also $0$, i.e., $\alpha^*=0$.
This indicates that there exist $t_i\in R, i=1, \dots,m, $ such that
 $\sum_{i=1}^m t_i A_i \succeq I \succ 0.$ ~  $\Box$

\vskip 0.08in

 As we mentioned in section 1, the result of Theorem
1.1 does not hold for $m\geq 3.$ That is, only knowing that the
system (\ref{LME-PSDC}) has no rank-one solution does not give a
full picture of the condition $\sum_{i=1}^m t_i A_i  \succ 0. $ The
system (\ref{LME-PSDC}) may have no rank-one solution, but have a
solution with rank$(X)\geq 2. $ For example, let {\small
\begin{equation} \label{EEE}  A_1= \left[
  \begin{array}{cccc}
    1 &   &   &      \\
     & -1 &   &     \\
      &   & 0 &     \\
     &  &  & 0    \\
  \end{array}
\right], ~ A_2=
  \left[
  \begin{array}{cccc}
  1 &   &   &      \\
     & 0 &   &      \\
      &   & -1 &     \\
     &  &  &  -1    \\
  \end{array}
\right], ~  A_3=
  \left[
  \begin{array}{cccc}
    0 &  1 &  0 &  0  \\
    1 & 0  &  0 &  0   \\
    0  &  0 & 0 & 0  \\
    0 & 0 & 0  & 0  \\
  \end{array}
\right].
\end{equation} }   It is easy to see that for this example
there is no $x\not=0$ satisfying $x^TA_ix=0, i=1,2,3,$  and hence
the corresponding system (\ref{LME-PSDC}) has no rank-one solution.
However, the system has a higher rank solution, $X=\textrm{diag}
(1,1,1,0).$ Clearly,  there exists no $(t_1,t_2, t_3) $ such that
$t_1A_1+t_2A_2+t_3A_3 \succ 0 $ for this example. To ensure the
condition $\sum_{i=1}^m t_i A_i \succ 0 ,$ Lemma 3.2 claims that the
system (\ref{LME-PSDC}) must possess not only no rank-one solution
but also
 no solution with rank higher than 1. It should be stressed that
the condition  $\sum_{i=1}^m t_i A_i \not\succ 0 $ for all $t_i$'s
implies that there is a nonzero solution $X$ to the system
(\ref{LME-PSDC}),
  but it cannot ensure that $
\textrm{rank}(X)=1 . $   Some stronger conditions than $\sum_{i=1}^m
t_i A_i \not\succ 0 $ should be imposed in order to guarantee the
existence of a rank-one solution.

From Lemma 3.2, we see that the linear combination of $A_i $'s plays
an important  role in determining the solution structure of the
 system (\ref{LME-PSDC}).  Given a finite number of matrices $A_i,
i=1,\dots,m,$ we use $r^*$ to denote the maximum rank of the linear
combination $\sum_{i=1}^m t_i A_i $, where $t_i\in R, i=1,\dots,m $
are chosen such that the linear combination is positive
semidefinite, i.e.,
\begin{equation}\label{R*}
r^*=   \max   \left\{  \textrm{rank}\left(\sum_{i=1}^m t_i
A_i\right) :  ~~  \sum_{i=1}^m t_i A_i\succeq 0, ~~ t_i \in R, ~
i=1,
 \dots,m  \right\}.
\end{equation}
  Clearly, $r^*$ is finite and attainable. Moreover, $r^*=n$ is equivalent
to $\sum_{i=1}^m t_i A_i \succ 0.$ The next result shows how  $r^*$
affects the existence of a nontrivial solution to the system
(\ref{LME-PSDC}), including low-rank ones.  It also indicates when
the rank-constrained  problem (\ref{rank-prob})  can be reduced to
an
SDP problem. \\

\textbf{Corollary 3.3.} (i) \emph{The system (\ref{LME-PSDC}) has a
solution $X\not=0$ if and only if  $r^*\leq n-1 $, where $r^*$ is
defined by (\ref{R*}). Moreover,  any nonzero solution  $X $ of
 (\ref{LME-PSDC}) satisfies $ {\rm rank}(X) \leq
n-r^*.$}

(ii)  \emph{Particularly, if $r^*=n-1$, then the  system
(\ref{LME-PSDC}) has a rank-one solution, and the rank-one solutions
are the only nonzero solutions of (\ref{LME-PSDC}). In this case,
the problem (\ref{rank-prob}) is equivalent to the SDP problem
\begin{equation}  \label {SDP-00}   \max   \{\langle I, X \rangle:
 ~\langle A_i, X \rangle = 0, ~i=1, \dots, m ,~
   \rm{tr}(X) \leq 1,~
  X\succeq 0\}.
  \end{equation}  }

  \vskip 0.08in

\emph{Proof. } (i) Lemma 3.2 claims that $r^*=n $ if and only if
$X=0$ is the only solution of  (\ref{LME-PSDC}).  Thus, the system
(\ref{LME-PSDC}) has a nonzero solution if and only if $r^*\leq n-1.
$  It is sufficient to prove that any nonzero solution $X$  of
(\ref{LME-PSDC}) must satisfy $ \textrm{rank} (X) \leq n-r^*.$
Indeed, let $(t_1^*,\dots, t_m^*) $ determine the maximum value
$r^*$, i.e.,
$$\label{rank}  r^* =\textrm{ rank} \left(\sum_{i=1}^m t_i^* A_i\right),
~~\sum_{i=1}^m t_i^* A_i \succeq 0.
$$
Let $X $ be an arbitrary solution of (\ref{LME-PSDC}). Thus, $X$
satisfies $$\left\langle X,~ \sum_{i=1}^m t_i^* A_i \right\rangle =
\sum_{i=1}^m t_i^* \left\langle X,~  A_i \right\rangle =0,
~~X\succeq 0. $$  Since $X\succeq 0 $ and $ \sum_{i=1}^m t_i^* A_i
\succeq 0,$ it implies that
\[
\left(\sum_{i=1}^m t_i^* A_i\right)X =0.
\]
Thus,  $ \textrm{rank}(X) \leq n-\textrm{rank}
\left(\sum_{i=1}^m t_i^* A_i\right) = n-r^*.$

We now prove (ii). Suppose $r^*=n-1$. Then the first half follows
 directly from (i) and Lemma 3.2.
Moreover, since any solution of (\ref{LME-PSDC}) satisfies
 $ \textrm{rank}(X)\leq n-r^*,$ it must satisfy
 $ \textrm{rank}(X)\leq 1. $ This means that the rank
constraint in (\ref{rank-prob}) is redundant. As a result, the
problem (\ref{rank-prob}) is reduced to the SDP problem
(\ref{SDP-00}).  ~~ $ \Box $

\vskip 0.08in

From Lemma 3.2 and Corollary 3.3, the cases $r^*=n$ and $r^* =n-1 $
are  clear. In the remainder of this section, we focus on the case
 $r^*\leq n-2 $ for which certain conditions should be imposed in order to ensure
the existence of a rank-one solution.
  In fact, when $r^*\leq n-2,$ the system
(\ref{LME-PSDC}) may have a solution with rank$(X) \geq 2 ,$ but no
rank-one solution.  It is easy to see that $r^*<n-2 $ holds in the
example (\ref{EEE}), since $r^*=0.$ Another simple example is that $A_1=
\left(
                    \begin{array}{cc}
                      1 & 0 \\
                      0 & -1 \\
                    \end{array}
                  \right) $ and $  A_2= \left(
                    \begin{array}{cc}
                      0 & 1 \\
                      1 & 0 \\
                    \end{array}
                  \right).$
For this example, we have $r^*=n-2$, but the system
 $\langle A_1, X \rangle=0,$ $\langle A_2, X \rangle=0,$ $X\succeq 0$ has
 no rank-one solution.
We now state an existence condition for the case $r^*\leq n-2. $

\vskip 0.08in

\textbf{Theorem 3.4.} \emph{Suppose that $r^*\leq n-2. $ Then the
system (\ref{LME-PSDC}) has a rank-one solution if any of the
following conditions holds:}

(i) \emph{For each $i$, $A_i$ is either positive semidefinite or negative semidefinite.}

(ii) \emph{There exists exactly one indefinite matrix among $A_i$'s.}

(iii) \emph{There are more than one indefinite matrices among
$A_i$'s, and there is an indefinite matrix $A_k$ such that
\begin{equation} \label{TTT} \{x: x^T
A_k x=0\} \subseteq \bigcap_{ {\rm All~indefinite}\,A_l, \, l\not= k}
\{x: x^T A_l x =0\}.
\end{equation} }

\emph{Proof.} When $r^*\leq n-2$, by Lemma 3.2 the system
(\ref{LME-PSDC}) has a solution with $1\leq \textrm{rank}(X) \leq
n-r^*.$ So if the system has no
  solution with rank~$(X) > 1$, then it must have a rank-one solution.
Thus, without loss of generality,
  we assume that the system (\ref{LME-PSDC}) has a solution $X^*$
with $\textrm{rank} (X^*) =r \geq 2, $ and hence $X^*$ can be decomposed as
$$X^*=\lambda_1 u^1(u^1)^T+  \cdots  +\lambda_r u^r(u^r)^T, $$
where $\lambda_j >0 $ and $u^j, j=1, \dots, r$
are eigenvalues and eigenvectors of $X^*, $ respectively, and
$u^j$'s are mutually orthogonal. For every $i=1,\dots, m, $ we have
$$ 0= \langle A_i, X^*\rangle =  \sum_{j=1}^r \lambda_j  (u^j)^T  A_iu^j  . $$
In particular, if $A_i \succeq 0$ or  $A_i \preceq 0$,
then the above equality implies
 \begin{equation}\label{LLL} (u^j)^T A_i
u^j = 0 ~ \textrm{ for all } j=1, \dots, r.
\end{equation}

First we suppose that condition (i) holds. That is, either  $ A_i
\succeq 0$ or  $A_i \preceq 0$ holds for each $i$. Then (\ref{LLL})
implies that any of the matrices $u^j (u^j)^T, j=1,\dots,r $ is a
rank-one solution to the system (\ref{LME-PSDC}).  This shows that
condition (i) ensures the existence of a rank-one solution.

Next, we suppose that condition (ii) holds and let $A_k$ be the only
indefinite matrix. If $ (u^j)^T A_k u^j =0$ for some $j$, then $X=
u^j(u^j)^T $ readily gives a rank-one solution to the system
(\ref{LME-PSDC}), since
 (\ref{LLL}) holds for all $i\not =k.$
On the other hand,  if $(u^j)^T A_k u^j \not =0  $ for all $j=1,
\dots, r $, then  from the fact that $\lambda_j >0$ for all $j$ and
$$
0= \langle A_k, X^*\rangle =  \sum_{j=1}^r \lambda_j (u^j)^T A_k u^j ,
$$ it follows that there exist two indices $p$ and $q$ such that
 $$ \left((u^p)^T A_k u^p\right) \left( (u^q)^T A_k u^q\right) < 0.$$
By continuity, there exists a $\gamma\in (0,1)$ such that
\begin{equation} \label{WWW1} w= \gamma u^p + (1-\gamma) u^q, ~~ w^T A_k w = 0 .
\end{equation}
Since $u^p$ and $u^q$ are orthogonal, it is evident that
$w\not=0. $  If $A_i \succeq 0 $ or  $A_i \preceq 0 $, then
(\ref{LLL}) implies that $A_i u^j=0$ for all $j=1, \dots, r.$
Thus we have
\begin{equation} \label{WWW2}
w^T A_iw = \gamma^2 (u^p)^T A_i u^p+(1-\gamma)^2 (u^q)^T A_i u^q
+2\gamma (1-\gamma) (u^p)^T A_i u^q =0 \end{equation} for all $i \ne
k$. Therefore, $X=ww^T $ is a rank-one solution of (\ref{LME-PSDC}).
Consequently, condition (ii) ensures the existence of a rank-one
solution.

Finally, we suppose that condition (iii) holds. If $ (u^j)^T A_k u^j
=0$ for some $j$, then it follows from (\ref{TTT}) that $(u^j)^T A_l
u^j =0$ for all other indefinite matrices $A_l.$ Since (\ref{LLL})
holds for all matrices $A_i$ such that $A_i\succeq 0 $ or
$A_i\preceq 0$, we may deduce that $X= u^j (u^j)^T $ is a rank-one
solution of the system (\ref{LME-PSDC}). If $(u^j)^T A_k u^j \not =0
$ for all $j=1, \dots, r$, then by the same reasoning as above, we
can find a vector $w \ne 0$ that satisfies (\ref{WWW1}). We can also
repeat the same argument as above to show that $w$ satisfies
(\ref{WWW2}) for all  $A_i\succeq 0 $ or $A_i\preceq 0$. Moreover,
by  (\ref{TTT}), we have $ w^T A_l w = 0 ~\textrm{ for all
indefinite } A_l \textrm{ with }  l\not = k. $ Thus, the nonzero
vector $w$ satisfies $w^T A_i w =0$ for all $i=1,\dots,m$, implying
that $X=ww^T $ is a rank-one solution to the system
(\ref{LME-PSDC}). The proof is complete. ~~ $ \Box $

\vskip 0.08in

 \textbf{Remark 3.5.} Given a set of matrices $A_i, i=1, ..., m,$
conditions (i) and (ii) in Theorem 3.4 can be verified straightaway.
A simple (and trivial) example satisfying the condition (iii) of
Theorem  3.4 is as follows:  Consider the system (\ref{LME-PSDC})
with $m=3$ and {\small $A_1=  \left(
          \begin{array}{ccc}
            1 & 0 & 0 \\
            0& 1 & 0 \\
             0& 0 & -1 \\
          \end{array}
        \right), $ } $     A_2 = 2 A_1 $ and $ A_3= 3 A_1. $ Then the
        condition (\ref{TTT}) holds trivially, and {\small  $X=
\left(
                                        \begin{array}{ccc}
                                          0 & 0 & 0 \\
                                          0 & 1 & 1 \\
                                          0 & 1 & 1 \\
                                        \end{array}
                                      \right) $ } is a rank-one
                                      solution of the system
                                      (\ref{LME-PSDC}).

\vskip 0.08in

 \textbf{Remark 3.6.} Theorem  3.4  shows that the number
of indefinite matrices among $A_i$'s and their relationships are
closely related to the existence of a rank-one solution to the
system  (\ref{LME-PSDC}). While this result gives some sufficient
conditions for the system (\ref{LME-PSDC}) to have a rank-one
solution, it is worth noting that these conditions remain not tight,
as shown by the following example: Consider the system
(\ref{LME-PSDC})  with {\small  $ A_1= \left(
          \begin{array}{ccc}
            1 & 0 & 0 \\
            0& 1 & 0 \\
             0& 0 & -1 \\
          \end{array}
        \right)$ }  and {\small $  A_2= \left(
          \begin{array}{ccc}
            1 & 0 & 0 \\
            0& -1 & 0 \\
             0& 0 & 1 \\
          \end{array}
        \right)$ }
which are both indefinite. It is easy to see that $$ \{x: x^TA_1 x
=0\} \not\subseteq \{x: x^T A_2 x=0\}, ~\{x: x^TA_2 x =0\}
\not\subseteq \{x: x^T A_1 x=0\}. $$ So all conditions (i), (ii) and
(iii) in Theorem 3.4 do not hold for this example. However, the
system (\ref{LME-PSDC}) with these two matrices has a rank-one
solution, for instance {\small $X= \left(
                                        \begin{array}{ccc}
                                          0 & 0 & 0 \\
                                          0 & 1 & 1 \\
                                          0 & 1 & 1 \\
                                        \end{array}
                                      \right) $ } is a rank-one
                                      solution.

\vskip 0.08in
 From Corollary 3.3(i), the rank of the solution of
(\ref{LME-PSDC}) is at most $n-r^*,$ which is a uniform bound for
all solutions. From a practical viewpoint, it is important to
compute the value $r^*.$ This motivates us to study the rank
maximization problem (\ref{R*}), which can be rewritten as
\begin{equation} \label{R**} r^* = \max\left\{ \lambda : ~
\textrm{rank} \left(\sum_{i=1}^m t_i A_i \right) \geq \lambda,
~\sum_{i=1}^m t_i A_i \succeq 0  \right\}.
\end{equation}
Since $\textrm{rank}(X)$ is a discontinuous function (in fact, a lower semi-continuous function)
of $X,$  the set $\{X: \textrm{rank}(X) \geq \lambda\}$ is not closed in general.
  This makes the problem (\ref{R**}) (or (\ref{R*}))
difficult to solve directly.  In what follows, we propose a method
to estimate $r^* $ from below. The following lemma will be used in
our analysis.

\vskip 0.08in

\textbf{Lemma 3.7.} \cite{F02, FHB04, RFP07} \emph{ The convex
envelope of $ {\rm rank} (X)$ on the set $\{ X\in R^{m\times n}:
\|X\|\leq 1\} $ is the nuclear norm $ \|X\|_*.$}

\vskip 0.08in

Note that for any matrix $Y\not =0 $, we have
$\textrm{rank}(Y)=\textrm{rank}(\alpha Y)$ for any $\alpha \not =
0.$ Thus problem (\ref{rank-prob}) can be rewritten as
 \begin{equation}\label{estimate} r^* = \max\left\{  \textrm{rank} \left(\sum_{i=1}^m t_i A_i
\right):  ~  \sum_{i=1}^m t_i A_i \succeq 0, ~ \sum_{i=1}^m t_i A_i
\preceq I \right\}.
\end{equation}
Since $ 0 \preceq \sum_{i=1}^m t_i A_i \preceq I$, we have
$\|\sum_{i=1}^m t_i A_i \|\leq 1 . $ By Lemma 3.7, we conclude that
in the feasible region of the problem (\ref{estimate}), the nuclear
norm of $\sum_{i=1}^m t_i A_i$  is the convex envelop  of the
objective function of (\ref{estimate}). As a result,  we have
\begin{equation} \label{nuclear}
 \textrm{rank}\left(\sum_{i=1}^m t_i A_i\right) \geq \left\|\sum_{i=1}^m t_i
A_i \right\|_*    \end{equation} for any $(t_1,\dots,t_m)$ satisfying $ 0
\preceq \sum_{i=1}^m t_i A_i \preceq I.$   By the positive
 semidefiniteness  of $ \sum_{i=1}^m t_i A_i,$ we have
$$ \left\|\sum_{i=1}^m t_i A_i \right\|_* =
\textrm{tr} \left(\sum_{i=1}^m t_i A_i \right) = \sum_{i=1}^m t_i \, \textrm{tr}(A_i).$$
Thus, we may consider the following  problem:
\begin{equation}\label{restric}  \eta^* = \max\left\{   \sum_{i=1}^m t_i \,
\textrm{tr}(A_i):  ~  \sum_{i=1}^m t_i A_i \succeq 0, ~ \sum_{i=1}^m
t_i A_i \preceq I \right\},
\end{equation}
which is an SDP problem with a finite optimal value $\eta^*\leq n.$
By (\ref{nuclear}) and Lemma 3.7, the optimal objective value of
(\ref{restric})  provides a  lower bound for that of
(\ref{estimate}), i.e., $r^* \geq \lceil\eta^* \rceil. $ The dual of
(\ref{restric}) is given by
\begin{equation} \label{D-SDP}
\min \{   \langle I, X \rangle:
  \langle A_i, X\rangle-\langle A_i,
Y\rangle=\textrm{tr}(A_i), ~ i=1,\dots, m,  ~X\succeq 0, ~Y\succeq
0\}.
 \end{equation}
This problem is strictly feasible and has a finite optimal value.
For instance, $(X,Y)=(2I,I)$ is a strictly feasible point. By the
duality theory, there is no duality gap between (\ref{restric}) and
(\ref{D-SDP}), and hence we may solve either of them to get the
optimal value $\eta^*.$ An immediate consequence of the above
analysis  is the following result. \\

\textbf{Theorem 3.8.} \emph{Let $r^*$ be the maximum rank defined by
(\ref{R*}), and let $\eta^* $ be the optimal value of the SDP
problem (\ref{restric}) or (\ref{D-SDP}).  When $r^*\leq n-2$, every
nonzero solution $X$ of the system (\ref{LME-PSDC}) satisfies $ {\rm
rank}(X)
\leq n- \lceil\eta^* \rceil. $} \\

This result provides an upper bound for the rank of nonzero
solutions of (\ref{LME-PSDC}), and the bound $n- \lceil\eta^* \rceil
$ can efficiently be obtained by solving  (\ref{restric}) or
(\ref{D-SDP}). \\

\textbf{Remark 3.9.} Consider the problem of finding $X\in S^n $
that satisfies
$$   \langle A_i, X \rangle = b_j, ~ j=1, \dots, m, ~ X \succeq 0. $$
If the above system has a solution, then it has a solution $X$ such
that
\begin{equation}  \label{BBBB} \textrm{rank}(X) \leq \left\lfloor \frac{\sqrt{8m+1}-1}{2}
\right\rfloor ,
\end{equation}
which is called Barvinok-Pataki's bound \cite{B95, PA98}. For a
homogeneous system (i.e., $b_j=0, j=1,\dots,m$), a solution
satisfying the bound (\ref{BBBB}) can be only the trivial solution
$X=0,$ and any nontrivial solution may not satisfy this bound. In
other words, \emph{the Barvinok-Pataki's bound is not necessarily
valid for nontrivial solutions of  a homogeneous system.} For
example, let $A_1= \left(
                    \begin{array}{cc}
                      1 & 0 \\
                      0 & -1 \\
                    \end{array}
                  \right) $ and $  A_2= \left(
                    \begin{array}{cc}
                      0 & 1 \\
                      1 & 0 \\
                    \end{array}
                  \right).
 $  The right-hand side of (\ref{BBBB}) is equal to 1 (since $m=2 $).
However,
  as we mentioned earlier, all nontrivial solutions to the system
 $\langle A_1, X \rangle=0, \, \langle A_2, X \rangle=0, \, X\succeq 0$ have
 rank~2. Thus, Barvinok-Pataki's  bound (\ref{BBBB})  cannot directly apply to
 nontrivial low-rank solutions of a homogeneous system
like  (\ref{LME-PSDC}). In order to apply this bound to the
homogenous system (\ref{LME-PSDC}), we may introduce an extra
equation, for instance,
 $\langle I, X \rangle=1,$ and consider the following system:
\begin{equation} \label{SLME}
\langle A_i, X \rangle =0, ~ j=1, \dots, m, ~\langle I, X \rangle=1,  ~ X \succeq 0.
\end{equation}
Note that $\textrm{rank} (tX)=\textrm{rank}(X)$ for any $t\not=0.$
By the homogeneity of (\ref{LME-PSDC}), we see that any nonzero
solution (if exists) of (\ref{LME-PSDC}) can be scaled so that it
satisfies (\ref{SLME}). Thus, if the system (\ref{LME-PSDC}) has a
nonzero solution, then the minimum rank of nonzero solutions of
(\ref{LME-PSDC}) and (\ref{SLME}) are the same. Therefore, applying
(\ref{BBBB}) to (\ref{SLME}), we can conclude from Theorem 3.6 that
when $r^*\leq n-2$, the system (\ref{LME-PSDC}) has a solution $X
\not= 0$ satisfying
\begin{eqnarray*}
1\leq \textrm{rank}(X) & \leq & \min \left\{n-r^*,
\left\lfloor \frac{\sqrt{8(m+1)+1}-1}{2}
\right\rfloor\right\} \\
&  \leq  & \min \left\{n- \lceil\eta^* \rceil, \left\lfloor
\frac{\sqrt{8(m+1)+1}-1}{2} \right\rfloor\right\}.
\end{eqnarray*} So when
$r^*$ is relatively large, e.g., $n-2, n-3, $ and so on, the rank of
a nontrivial solution to (\ref{LME-PSDC}) will be low.  In such
cases,  Barvinok-Pataki's bound might be too loose (especially when
$m$ is relatively large). The
 upper bound given by $n-r^*$ or even $n- \lceil\eta^* \rceil$  for
the rank of nonzero solutions can be much tighter than Barvinok-Pataki's bound
in these situations.

\section{Conditions for (\ref{LME-PSDC}) to have no rank-one solution}

The following necessary condition for the system (\ref{LME-PSDC}) to
have no rank-one solution has actually been shown in section 2 (see
Lemmas 2.1 and 2.4, or Corollary 2.5).

\vskip 0.08in

\textbf{Corollary 4.1.} \emph{When $m\leq n-1,$ if the system
(\ref{LME-PSDC}) does not have a rank-one solution, then there
exists a  vector $x$ such that $$ x^T A_i x= \frac{x_i}{x_n}, ~
i=1,\dots, m , ~ \|x\| =1, ~   x_n\in (-1,0).
$$}

 Although the condition $\sum_{i=1}^m t_i A_i \succ 0$
ensures that the system (\ref{LME-PSDC}) has no rank-one solution,
this sufficient condition is too restrictive. In fact, by Lemma 3.2,
it implies that the system (\ref{LME-PSDC}) cannot have any solution
with rank$(X)\geq 1.$
 The  purpose of this section is to show that another sufficient condition
 for  the system (\ref{LME-PSDC}) to have no rank-one solution can be developed from
 a homotopy invariance point of view.
Note that the following three statements are equivalent: (i) The
system (\ref{LME-PSDC}) has no rank-one solution; (ii)  $x=0$ is the
only solution to the system (\ref{quadratic}); and (iii)
 \begin{equation}\label{positive}
\max_{1\leq i\leq m}  |x^TA_i x| >0 ~~ \textrm{for any } x\not=0.\end{equation}
 First, we formulate these equivalent statements as a nonlinear
 equation.

 \vskip 0.08in

\textbf{Lemma 4.2.}  \emph{$x=0$ is the only solution to the system
(\ref{quadratic}) if and only if there exists a constant $\beta>0 $
 such that for any $\mu \in (0, \beta],$ we have
$$ \{x: \|x\|=1\} = \{x: G_\mu(x)=0 \},$$ where $G_\mu: R^n \to R^2 $ is defined by
\begin{equation} \label{GGG-111}
G_\mu (x)  = \left(\begin{array}{cc}  \left| \sum_{i=1}^m |x^T A_i x|-\mu \right|
             -\left(\sum_{i=1}^m |x^T A_i x|-\mu  \right) \\ [3pt]
x^T x-1
\end{array}\right), \nonumber  \end{equation}
i.e.,  the set $\{x: \|x\|=1\}$ coincides with the solution set of
the equation $G_\mu(x) =0 $ for any $\mu\in (0, \beta].$}

\vskip 0.08in

\emph{Proof.}  Assume that $x=0$ is the only solution of the system
(\ref{quadratic}). Thus, by (\ref{positive}),  we have $\sum_{i=1}^m
|x^T A_i x| > 0 $ for any $x$ such that $x^Tx=1.$  By continuity,
there exists a positive number $\beta>0$ (for instance, we can take
$\beta= \min \left\{
 \sum_{i=1}^m |x^T A_i x| : ~ \|x\|=1\right\} $ which is
positive) such that
\begin{equation} \label{EEE-01} \sum_{i=1}^m |x^T A_i x| \geq \mu
  \end{equation} holds for any    $ x $    with  $  x^T x=1$  and
for any $\mu \in (0, \beta].$
Note that any inequality $h(x)\geq 0$ can be represented as the
equation $|h(x)|-h(x) =0.$ So (\ref{EEE-01}) can be rewritten as
$$\left| \sum_{i=1}^m |x^T A_i x|-\mu \right| -
\left(\sum_{i=1}^m |x^T A_i x|-\mu  \right) =0.   $$
This implies that any $x$ satisfying $x^Tx=1 $ is a solution to the
equation
 $$G_\mu(x) =0 $$
for any $\mu \in (0,\beta]$.
Thus, the set $\{ x : \|x\|=1 \}$ is contained in the solution set of
$G_\mu(x)=0 $ for any  $\mu \in (0,\beta]$. Since $G_\mu(x) =0 $
implies $\|x\|=1,$   the set $\{ x : \|x\|=1 \}$ is exactly the solution set
of $G_\mu(x)=0 $  for any  given $\mu \in(0, \beta].$

Conversely, let us assume that there exists a positive number $\beta
>0$ such that for any given constant $\mu \in (0, \beta]$, the
solution set of the equation $G_\mu(x)=0 $ is equal to the set $\{ x
: \|x\|=1 \}$. We now prove that $x=0$ is the only solution of
(\ref{quadratic}). Assume the contrary that the system
(\ref{quadratic}) has a solution $x\not= 0. $  Then,
$\widehat{x}=x/\|x\| $ is also a solution of the system
(\ref{quadratic}), i.e., $\widehat{x}^T A_i \widehat{x} =0, ~i=1,
\dots,m. $ By assumption, any unit vector is a solution to
$G_\mu(x)=0 $ for any given $\mu \in (0, \beta].$ Thus, we have
$$ \left| \sum_{i=1}^m |\widehat{x}^T A_i \widehat{x}|
 -\mu  \right| -\left(\sum_{i=1}^m |\widehat{x}^T A_i
\widehat{x}|-\mu     \right) =0. $$ Since
$\widehat{x}^T A_i \widehat{x} =0$ for all $ i=1, \dots, m$, the above
equality reduces to
$$  0 =2\mu  , $$
which is a contradiction since $  \mu \in(0, \beta]. $  ~  $ \Box$
\\

We assume $m\leq n-2$ in the remainder of this section.
  Again, by using Lemma 2.2, we have the following technical result.
  \\

\textbf{Lemma 4.3. } \emph{Assume that $m\leq n-2.$ Let $\mu>0 $ be a
given constant and $G_\mu $ be defined by (\ref{GGG-111}). If
$G_\mu(x^*) \not = 0 $ for some $x^*$ with $\|x^*\| =1 , $ then
  there exists a sequence $\{x^k\}  $ satisfying the
  following conditions:  $\|x^k\|\to \infty  $ as $k\to \infty$ and
\begin{eqnarray} (x^k)^T A_1 x^k - (x^*)^TA_1 x^* & = & -\gamma^k x^k_1,  \label{first}\\
 & \vdots & \nonumber\\
 (x^k)^T A_m x^k  - (x^*)^T A_m x^* & = &  -\gamma^k x^k_m , \label{second}
\end{eqnarray}
 \begin{equation} \label{last-second}
 \left| \sum_{i=1}^m |(x^k)^T A_i x^k| -\mu  \right|
- \left(\sum_{i=1}^m |(x^k)^T A_i x^k|-\mu \right) =  - \gamma^k x^k_{n-1} ,
 \end{equation}
 \begin{equation} \label{last} (x^k-x^*)^T (x^k-x^*)     =    -\gamma^k x^k_{n},
 \end{equation}
where $\gamma^k\in (0, \infty) $ for all $k.$}

 \emph{Proof.} If $m< n-2$, we may set $A_{m+1}=0, \dots, A_{n-2}=0,$ and consider the systems
 (\ref{LME-PSDC}) and (\ref{quadratic}) with $A_i, i=1, \dots, m \, (=n-2).$
Thus, without loss of generality, we assume $m=n-2 $.
  Since $G_\mu(x^*)\not =0 $ and $\|x^*\|=1,  $  we have  $$ \left| \sum_{i=1}^m |(x^*)^T A_i
x^*|-\mu   \right| -\left(\sum_{i=1}^m |(x^*)^T A_i
x^*|-\mu   \right) \not= 0.$$ Thus the following
equation has no solution:
\begin{equation} \label{GGG-333}
\Theta_\mu (x) = \left(\begin{array}{cc} x^TA_1 x-(x^*)^TA_1 x^* \\
 \vdots\\
  x^TA_m x-(x^*)^TA_m x^*\\  \left|
\sum_{i=1}^m |x^T A_i x|-\mu    \right| -\left(\sum_{i=1}^m
|x^T A_i x|-\mu
 \right) \\
(x-x^*)^T (x-x^*)\\
\end{array}\right) =0.  \end{equation}
Consider the homotopy between the identity mapping and $\Theta_\mu
(x) $, i.e.,
$$ H_\mu(x,t)= t x +(1-t)\Theta_\mu (x),
$$
and let
$$ {\mathcal T_\mu}  = \{x \in R^{n}: ~ H _\mu(x,t)=0\textrm{  for some  } t \in [0,1] \} .  $$
A similar proof to that of Lemma 2.3 can be used to show that
${\mathcal T_\mu}$ is unbounded. Here we include the proof for
completeness. Assume to the contrary that ${\cal T}_\mu $ is bounded.
Then there exists an open bounded ball $D$ that is large enough to satisfy
${\mathcal T}_\mu \subset D $ and
$\partial D\bigcap \mathcal{T_\mu} =\emptyset.$
Thus $0 \notin\{H_\mu (x,t): ~x\in
\partial D, ~0\leq t\leq 1\},$ which implies that $\deg
(I, D,0) $ and $\deg(\Theta_\mu , D, 0)$ are well defined. By Lemma
2.2, we have
$$
\deg(\Theta_\mu , D, 0)=  \deg (I, D,0)  \not = 0,  $$ which means
that the equation $ \Theta_\mu (x)=0  $  has a solution. This is a
contradiction. Thus the set $\mathcal{T_\mu} $ is indeed unbounded,
and hence there is an unbounded sequence $\{x^k\}$
in $\mathcal{T}_\mu$. Without loss of generality,
let $x^k\ne 0$ for all $k.$ By the definition of $\mathcal{T}_\mu$,
there is a sequence $\{t^k\}\subseteq [0,1] $ such that $H_\mu(x^k,
t^k) =0$, i.e.,
 \begin{eqnarray*}
&  t^kx^k_1+ (1-t^k) [(x^k)^T A_1 x^k-(x^*)^TA_1 x^*]  =0, &  \\
&     \vdots  & \\
& t^kx^k_m+ (1-t^k) [(x^k)^T A_m x^k-(x^*)^TA_m x^*] =0, &
\end{eqnarray*}
\begin{eqnarray*} & &   t^k x^k_{n-1} +
  (1-t^k)
  \left[\left| \sum_{i=1}^m |(x^k)^T A_i
x^k|-\mu   \right|
    -\sum_{i=1}^m |(x^k)^T A_i x^k|
       +\mu
 \right] =0,
 \end{eqnarray*}
 $$   t^k x^k_{n} +(1-t^k)  (x^k-x^*)^T (x^k-x^*)  =0 .$$
Since $\|x^k\|\to \infty $ as $k\to \infty$ and  $\Theta_\mu (x^k)
\not=0 $ for all $k\geq 1,$ we see that $ t^k\not=1 $  and
$t^k\not=0 $ for all $k\geq 1. $   Thus, $t^k\in (0,1)$ for all $k
\geq 1. $ By setting $\gamma^k = \frac{t^k}{1-t^k} \in (0, \infty) ,
$  the above system can be rewritten as
$$ \Theta_\mu  (x^k) = -  \frac{t^k}{1-t^k}x^k = -\gamma^k  x^k, ~~ k\geq 1, $$
which along with the definition of $ \Theta_\mu $ implies that
\begin{eqnarray*} (x^k)^T A_1 x^k  -(x^*)^TA_1 x^*& = & -\gamma^k x^k_1,  \\
& \vdots & \\
 (x^k)^T A_m x^k -( x^*)^T A_m x^* & = & -\gamma^k x^k_m,
 \end{eqnarray*}
 $$   \left| \sum_{i=1}^m |(x^k)^T A_i
x^k|-\mu   \right| -\left(\sum_{i=1}^m |(x^k)^T A_i
x^k|-\mu
 \right)  =    - \gamma^k x^k_{n-1}, $$
 $$ (x^k-x^*)^T (x^k-x^*)     =     -\gamma^k x^k_{n},$$
as desired.  $ ~  \Box $

\vskip 0.08in

Based on Lemma 4.3, we can prove the following result.

\vskip 0.08in

\textbf{ Lemma 4.4.} \emph{Let  $m\leq n-2$ and $\mu>0$ be a given
constant. If $G_\mu(x^*) \not=0$ for some $x^*$ with $\|x^*\|=1$,
then there exists a vector $ \widehat{x} $  satisfying the
following conditions:
\begin{equation}\label{condition_in_Lemma4.4}
  \widehat{x}^T A_i \widehat{x}    =   \frac{\widehat{x}_i }{ \widehat{x}_{n} },
 ~ i=1, \dots, m , ~ \widehat{x}_{n-1}=0 ~ \widehat{x}_{n} \in [-1,0),
~ \|\widehat{x}\|=1.
\end{equation} }

\vskip 0.08in

 \emph{Proof.}  Without loss of generality, we still assume that $m=n-2.$
Let $\{x^k\}$ be the sequence specified in Lemma 4.3.
Since the left-hand sides of (\ref{last-second})
  and (\ref{last}) are nonnegative, we see that
 $x^k_n <0 $ for all $k\geq 1.$
  Let $\widehat{x}$ be an accumulation point of $x^k/\|x^k\|.$
By dividing (\ref{first})--(\ref{last}) by $\|x^k\|^2$ and
passing to a subsequence if necessary,  we can prove  that $
x^k_n/\|x^k\|\to \widehat{x}_{n} \not =
 0.$
In fact, if $\widehat{x}_{n}  = 0$,
 (\ref{last}) implies $ \gamma^k/\|x^k\| \to \infty. $
Then it follows from (\ref{first})--(\ref{last-second}) that
 $\widehat{x}_i=0$ for all $ i=1,\dots,n-1.$ This contradicts the fact
 $\|\widehat{x}\|=1.$
 Therefore, we conclude that $\widehat{x}_{n} \in
 [-1,0).$ It then follows from (\ref{last}) that $ \gamma^k/\|x^k\|
 \to \widehat{\gamma} \in (0, \infty), $  where $\widehat{\gamma}=-1/\widehat{x}_n$.
Normalizing the
 system (\ref{first})--(\ref{last}) by $\|x^k\|^2 $ and letting $k\to \infty$
 yield
$$ \widehat{x}^T A_i \widehat{x}     =
 -\widehat{\gamma}  \widehat{x}_i , ~ i=1, \dots, m ,
 ~0 =   -\widehat{\gamma} \widehat{x}_{n-1},
~1 = -\widehat{\gamma}  \widehat{x}_{n}    ,
 ~ \widehat{x}_{n}   \in   [-1,0).
$$
Eliminating $\hat{\gamma}$ from the above system yields the desired
result.  ~~ $  \Box $

\vskip 0.08in

Basically, the next result shows that what an extra condition can
make the necessary condition in Corollary 4.1 sufficient.

 \vskip 0.08in

\textbf{Theorem 4.5.} \emph{Let $m \leq n-2 $ and $A_i \in S^n,
~i=1,\dots, m. $ Suppose that $e_n=(0,\dots,0,1)^T\in R^n$ is not a
solution of the system (\ref{quadratic}), and suppose that there is
a vector $x$ satisfying the following condition:
\begin{equation}\label{sufficient}
\label{infeasible} x^T A_i x = \frac{x_i}{x_n} , ~ i=1, \dots, m,  ~
x_{n} \in (-1, 0), ~ \|x\| =1. \end{equation} If $x_{n-1} \not =0 $
for any $x$ satisfying (\ref{sufficient}),
 then $x=0$ is
the only solution to the system (\ref{quadratic}), i.e., the system
(\ref{LME-PSDC})
 has no rank-one solution.}

 \vskip 0.08in

\emph{Proof. } Under the condition of this theorem, we prove that
there is a $\beta>0$ such that
$$\{x:\|x\|=1\} = \{x: G_\mu(x)=0\} \quad
\textrm{ for any } \, \mu \in (0, \beta],$$ where $G_\mu$ is defined
by (\ref{GGG-111}), and thus by Lemma 4.2,  $x=0$ is the only
solution to the system (\ref{quadratic}). We prove this  by
contradiction. Assume that such a $\beta$ does not exist, i.e., for
any given $\beta >0 $ (no matter how small it is),  there always
exists a $\mu \in (0, \beta]$ such that   $\{x:\|x\|=1\} \not = \{x:
G_\mu(x)=0\}.$ In other words, there exists an $x^*$ with $\|x^*\|=1
$ such that $G_\mu(x^*) \not =0. $ Then, by Lemma 4.4, we conclude
that
 there exists an $\widehat{x}$ satisfying (\ref{condition_in_Lemma4.4}).
Since  $e_n=(0,...,0,1)^T\in R^n$ is not a solution of the system
(\ref{quadratic}), we deduce that $\widehat{x}_{n}\not= -1$ in
 (\ref{condition_in_Lemma4.4}), and hence $\widehat{x}$
satisfies (\ref{sufficient}) and $\widehat{x}_{n-1}=0,$
  contradicting the assumption of the theorem.  ~~ $   \Box $ \\

\textbf{Remark 4.6.} The above theorem provides a sufficient
condition for the system (\ref{LME-PSDC}) to have no rank-one
solution. This is equivalent to saying that when the system
(\ref{LME-PSDC}) has a rank-one solution, such a sufficient
condition must fail. For example, let us consider the following
system: {\small
\begin{equation} \label{LAST-EX}  A_1= \left[
  \begin{array}{ccccc}
    1 &   &   &   &   \\
     & 0 &   &   &  \\
      &   & 0 &   & \\
     &  &  & -1  &   \\
    &  &  &   &  0
  \end{array}
\right], ~ A_2=
  \left[
  \begin{array}{ccccc}
    1 &   &   &   &   \\
     & 0 &   &   &  \\
      &   & -1 &    & \\
     &  &  & -1  &   \\
    &  &  &   &  -1
  \end{array}
\right], ~  A_3=
  \left[
  \begin{array}{ccccc}
    0 &  1 &  0 &  0  &  0\\
    1 & 0  &  0 &  0  &  0 \\
    0  &  0 & 0 & 0 & 0 \\
    0 & 0 & 0  & 0  & 0\\
 0 & 0 & 0  & 0  &  0\\
  \end{array}
\right].
\end{equation} }
For this example,  $e_5$ is not a solution of the system
(\ref{LME-PSDC}), and the condition (\ref{sufficient}) can be
written as
$$\left\{\begin{array}{l}
 x_5(x_1^2-x_4^2)=x_1 \\
 x_5(x_1^2-x_3^2-x_4^2-x_5^2) =x_2 \\
  x_5(x_1x_2)=x_3 \\
 x_5\in (-1,0), ~  \|x\|=1.
\end{array}\right.
$$
If we set $x_4=0$, then the above condition imply that $x_1=x_3=0,$
$x_2=-x_5^3$ and $ x_2^2+x_5^2=1.$ Therefore,   the point  $ x= (0,
-t^3, 0, 0,  t)^T$ satisfies the above condition with $$t= - \left(
\left( \frac{1}{2} +\sqrt{\frac{1}{4}+\frac{1}{27}}  \right)^{1/3} +
\left(\frac{1}{2} -\sqrt{\frac{1}{4}+\frac{1}{27}} \right)^{1/3}
\right)^{1/2}  \in (-1, 0) . $$ So for this example  the sufficient
condition in Theorem 4.5 does not hold, and it is easy to see that
the system  (\ref{LME-PSDC}) with matrices given by (\ref{LAST-EX})
has a rank-one solution, for instance, $X=xx^T$ where $x=
(0,1,0,0,0)^T.$

\vskip 0.08in

 \textbf{ Remark 4.7.}  From the analysis in this
paper, the condition (\ref{sufficient}) is intrinsically hidden
behind the condition ``there is no rank-one solution to the system
(\ref{LME-PSDC})", or equivalently, ``$x=0$ is the only solution to
the system (\ref{quadratic})".
 Theorem 4.5 shows that this type of necessary
 condition together with some other conditions can be sufficient for
 the nonexistence of a rank-one solution to the system (\ref{LME-PSDC}).  However, the
 relationship between the new sufficient conditions in this paper and
 the known condition ``$\sum_{i=1}^m t_i A_i\succ 0$" is not clear at
 present. These two types of conditions seem independent to each other.
Corollary 3.3 and Theorem 3.4 indicate that
 only a small gap exists between the existence and nonexistence of
a rank-one solution to the system (\ref{LME-PSDC}).
 If $r^* =n$, then $X=0$ is the only solution of the system
(\ref{LME-PSDC}). However, a small perturbation of the system  such
that $r^*=n-1$ will guarantee the existence of a rank-one solution
to the system.  Thus, the development of  a new sufficient condition
weaker than $\sum_{i=1}^m t_i A_i\succ 0$ becomes subtle, and there
might be no easy and simple way to state such a sufficient
condition. As we have shown in section 3, checking whether or not
the system (\ref{LME-PSDC}) has a rank-one solution is equivalent to
solving an SDP problem with rank constraints which is a difficult
problem. This  indicates that the conditions developed in sections 2
and 4 of this paper are not easy to check directly. However, these
conditions make it possible for us to understand the problem from
the nonlinear analysis perspective.

\section{Conclusion}
Some sufficient and/or necessary conditions for the existence of a
rank-one solution to the system of HLME over the positive
semidefinite cone have been developed. These conditions have been
derived from two different perspectives: degree theory  and  rank
optimization. The result out of the former shows that the
$P_\emptyset$ property of the function defined by quadratic
transformations can ensure the existence of a rank-one solution
(e.g., Theorems 2.6 and 2.9). From the latter, it turns out that the
maximum rank $r^*$, defined by (\ref{R*}), plays a key role in the
existence of a rank-one solution. For instance, $r^*=n-1$ can ensure
the system of HLME has a rank-one solution (see Corollary 3.3), and
the number of indefinite matrices in the system can  be also related
to the existence of a rank-one solution (see Theorem 3.4).   Finally
a sufficient condition for the nonexistence of a rank-one solution
to the system of HLME was also given (see Theorem 4.5).

\end{document}